\documentclass[CJK,11pt]{article}
\usepackage{amscd}
\usepackage{mathrsfs}
\usepackage{amsfonts}
\usepackage{amsmath,amscd,stmaryrd,amssymb,array, amsfonts,amsmath,amssymb,graphicx}

\numberwithin{figure}{section}
 \numberwithin{equation}{section}

\newtheorem{theorem}{Theorem}[section]
\newtheorem{proposition}[theorem]{Proposition}
\newtheorem{definition}[theorem]{Definition}

\newtheorem{lemma}[theorem]{Lemma}
\newtheorem{remark}[theorem]{Remark}

\newcommand{\cA}{{\mathcal A}}
\newcommand{\cB}{\mb{B}}

\newcommand{\cC}{{\mathcal C}}
\newcommand{\cD}{{\mathcal D}}

\newcommand{\cG}{{\mathcal G}}

\newcommand{\cM}{{\mathcal M}}

\newcommand{\cV}{{\mathcal V}}
\newcommand{\cR}{{\mathcal R}}
\newcommand{\cS}{{\mathcal S}}

\newcommand{\cK}{{\mathcal K}}

\def\be{\begin{equation}}
\def\ee{\end{equation}}
\def\ba{\begin{array}}
\def\ea{\end{array}}
\def\benu{\begin{enumerate}}
\def\eenu{\end{enumerate}}
\def\bt{\begin{theorem}}
\def\et{\end{theorem}}
\def\bp{\begin{proposition}}
\def\ep{\end{proposition}}

\def\bl{\begin{lemma}}
\def\el{\end{lemma}}
\def\br{\begin{remark}}
\def\er{\end{remark}}

\def\b{\beta}
 \def\de{\delta} \def\pa{\partial} \def\nab{\nabla}
 \def\O{\mathcal{O}}
\def\ve{\varepsilon} \def\sig{\sigma}
\def\w{\omega}\def\W{\Omega}
\def\a{{\mathscr A}}

\def\.{\cdot}
\def\cdt{\stackrel{\centerdot}{}}

\def\R{\mathbb{R}}

\def\A{\forall}
\def\ol{\overline}

\def\Cap{\bigcap}\def\Cup{\bigcup}
\def\ln{\mbox{ln\,}}
\def\ra{\rightarrow}
\def\stac{\stackrel}
\def\~{\stac{\sim}}
\def\8{\infty}
\def\X{\times}
\def\({\left(}
\def\){\right)}
\def\d{\partial}

\def\mb{\mbox}
\def\emp{\emptyset}
\def\-{\setminus}

\def\Hs{\hspace{1cm}}\def\hs{\hspace{0.5cm}}
\def\Vs{\vskip8pt}\def\vs{\vskip4pt}

\def\({\left(}\def\){\right)}

\begin{document}

\begin{center}{\Large\bf Smooth Morse-Lyapunov Functions and\\
Morse Theory of Strong Attractors for\\[1ex] Nonsmooth Dynamical Systems
 }\footnote{This work is supported by
NNSF of China (10771159) and NSF of Tianjin (09JCYBJC01800).
\vskip2pt \hspace{8pt}2010 {\em  Mathematics Subject
Classification.} 37B25, 37B30. \vskip2pt \hspace{8pt}{\em Key words
and phrases.} Nonsmooth dynamical system, attractor, Morse-Lyapunov
function, Morse theory. }

 \end{center}
 \centerline{\scshape  Desheng Li} \vskip4pt { \centerline{Department of
Math., Tianjin University} \centerline{ Tianjin 300072, P.R.
 China}
  \centerline{\footnotesize {\em Email}:\, lidsmath@tju.edu.cn, lidsmath@hotmail.com} }
\vskip10pt \centerline{{\em Date\,}: May 18, 2010}

\vskip10pt

\begin{center}
\begin{minipage}{12cm}
\Vs{\bf Abstract.} In this paper we first construct smooth
Morse-Lyapunov functions of attractors for nonsmooth dynamical
systems. Then we prove that all open attractor neighborhoods of an
attractor have the same homotopy type. Based on this basic fact  we
finally introduce the concept of critical group for Morse sets of
 an attractor  and establish Morse inequalities and equations. \Vs

\Vs \Vs \centerline{\bf Contents} \Vs
\begin{enumerate}

\item[1.] Introduction

\item[2.] Preliminaries

\item[3.] Smooth Lyapunov Functions of Attractors

\item[4.] Smooth Morse-Lyapunov Functions

\item[5.]  Morse Theory
\end{enumerate}
\end{minipage}
\end{center}

\newpage

\section{Introduction}\label{intro}

Attractors are of particular interest in the theory of dynamical
systems, this is because that much of the longtime dynamics of a
given system  is represented by the dynamics on and near the
attractors.

 The attractors of smooth
dynamical systems have been extensively studied in the past decades,
in both finite and infinite dimensional cases. The existence results
(especially  for infinite dimensional systems) are well known; see
 \cite{CV,Lady,Rob,Sell} and \cite{Tem} etc.
 In many cases one can
 give an estimate on the
  Hausdorff (or fractional, or informational)
dimension of an attractor. It can even be proved that the global
attractor of an infinite dimensional system is actually contained in
a finite-dimensional manifold; see \cite{Const,Tem} for details. The
Morse theory  is also fully developed \cite{Conley,CZ,Kap,
Kur,Ras,RS,Ryba,San}.
 In contrast, the situation in nonsmooth dynamical
systems seems to be more complicated, and by far some fundamental
problems concerning attractors are still undergoing investigations
\cite{CLV, CKM, Che, KV,Lisiam, mv,mv2}.

Nonsmooth dynamical systems appear widely in a large variety of
applications such as mechanics with dry friction, electric circuits
with small inductivity, systems with small inertial, economy,
biology, viability theory, control theory, game theory and
optimization etc. \cite{A,Aubin,AF, Bress2, Brog,CLS,Danca,
Deim,Feck1,Feck2, Fili,Glock, Hu,Lei,TP}. The rapid growth of such
systems in recent years challenges mathematicians to develop more
direct and uniform approaches to study their dynamics. In this paper
we are basically interested in finite dimensional case, in which
nonsmooth systems can be typically described  in a quite uniform
manner via  the following differential inclusion:
\begin{equation}\label{e:1.1}
x'(t)\in F(x(t)),\Hs x(t)\in X:=\R^m.
\end{equation}
One of the main feature of (\ref{e:1.1}) is that it may fail to have
uniqueness on solutions. Because of this, it usually generates a
multi-valued semiflow. So one needs to distinguish dynamical
concepts between strong and weak sense, where the former means that
they apply to ``all\,'' solutions, and in the latter ``all\,'' is
replaced by ``some\,''. The two settings are rather different. We
make precise that in this present work we will be solely interested
in the {\em strong} case. Hence from now on all the dynamical
concepts concerning (\ref{e:1.1}) should be understood in the strong
sense, except otherwise statement.

In our previous work  \cite{Lisiam} we have discussed Morse
decompositions of attractors for (\ref{e:1.1}). Morse decompositions
reveal some topological structures of  attractors, and are of
crucial importance in the understanding of the dynamics inside
attractors. Here  we want go a further step. First we construct
smooth Morse-Lyapunov functions for attractors. Then we prove that
all open attractor neighborhoods of an attractor have the same
homotopy type. Based on this basic fact we further introduce the
concept of critical groups for Morse sets and establish Morse
inequalities and equations of attractors.

Now let us give a brief description of what we will do. We will
assume and {\em only} assume  throughout the paper that $F$
satisfies the following {\bf standing assumptions}:

\vskip5pt
 (H1)
{ $F(x)$ is a nonempty convex compact subset of $X$ for every $x\in
X$;} \vs (H2) { $F(x)$ is upper semi-continuous in $x$.}

\vskip5pt

\noindent Let $\a$ be an attractor of (\ref{e:1.1}) with attraction
basin $\W=\W(\a)$ and Morse decomposition $\cM=\{M_1,\cdots,M_l\}$.
We will construct a radially unbounded smooth function $V\in
C^\8(\W)$ such that \Vs(1)
 $V$ is constant on each Morse set $M_k$, and\vs(2)
$V$ is strictly decreasing along  solutions of (\ref{e:1.1}) in $\W$
outside the Morse sets. \Vs \noindent Moreover, there is a
nonnegative function $w\in C(\W)$ which is positive on $\W\-\cD$,
where $\mathcal{D}=\cup_{1\leq k\leq l}M_k$, such that
\be\label{e:1.3} \max_{v\in F(x)} \nab V(x)\stackrel{\centerdot}{}
v\leq -w(x),\Hs \A\,x\in \W\-\mathcal{D}.\ee

Lyapunov functional characterizations of dynamical behavior are
usually known as converse Lyapunov theorems. This question can be
traced back to Lyapunov \cite{Lyap}. One of the early important
milestones in the pursuit of smooth converse Lyapunov functions was
Massera's 1949 paper \cite{Massera}. Since then a vast body of
fundamental results and techniques are obtained for smooth dynamical
systems (and even for ODE systems with only continuous righthand
sides \cite{Kurz}); see references cited in \cite{CL,TP} etc.
Nevertheless, the problem of smooth Lyapunov characterization of
stability for nonsmooth systems had remained open for a long time,
although Krasovskii pointed out as early as 1959 the desirability of
such characterizations \cite{Kras}.  A great progress in nonsmooth
case  was first made by Clarke, Ledyaev and Stern in their 1998
paper \cite{CL}, in which the authors proved a smooth converse
Lyapunov theorem for  strong global asymptotic stability of the zero
solution of (\ref{e:1.1}). Later on a more general framework
appeared in Tell and Praly's work \cite{TP}. Other related works can
 be found in \cite{LSW, ST} etc.

Note that an attractor $\a$ always has a trivial Morse decomposition
$\cM=\{\a\}$. In such a case a Morse-Lyapunov function of $\cM$
reduces to a  Lyapunov one  of the attractor itself. In particular,
in case $\a$ is an equilibrium we recover the converse Lyapunov
theorem in \cite{CL}, except that we don't require the function $w$
in (\ref{e:1.3}) to be smooth.
 To the best of
our knowledge, the above smooth converse Lyapunov theorem for Morse
decompositions is rather new even if we come back to the situation
of smooth dynamical systems.

We then try to  establish a Morse theory of attractors for nonsmooth
dynamical systems by employing smooth Lyapunov functions and
Morse-Lyapunov functions. In the smooth case this can be done by
using Conley index or shape theory
\cite{Conley,CZ,Kap,Kur,RS,Ryba,RZ}. Here we want to develop a
 somewhat different  approach for nonsmooth  systems. First, we prove that all the
open attractor neighborhoods of an attractor have the same homotopy
type. Then based on this basic fact we introduce the concept of
critical group for Morse sets and establish Morse equations and
inequalities. Specifically, let $$ \emp=A_0\subset
A_1\subset\cdots\subset A_l=\a$$ be the Morse filtration of the
Morse decomposition $\cM=\{M_1,\cdots, M_l\}$.  We define the
critical group $C_*(M_k)$ of a Morse set $M_k$ to be the homologies
of the space pair $(W,U)$ (with coefficients in a given Abelian
group $\mathcal{G}$) for any open attractor neighborhoods $W$ of
$A_k$ and $U$ of $A_{k-1}$, that is,
$$
C_*(M_k)=H_*(W,U).
$$
This definition is independent of the choice of $W$ and $U$. Set
$$
\mathfrak{m}_q=\sum_{1\leq k\leq l}\mb{rank}\,C_q\(M_k\),\Hs
q=0,1,\cdots.
$$
Let  $\b_q:=\b_q(\W)=\mb{rank}\,H_q(\W)$ be  the $q$-th Betti number
of the attraction basin $\W.$ We prove that the following Morse
inequalities and equation hold:
$$
\mathfrak{m}_0\geq \beta_0,$$
$$
\mathfrak{m}_1-\mathfrak{m}_0\geq \beta_1-\beta_0,
$$
$$
\cdots\cdots
$$
$$
\mathfrak{m}_m-\mathfrak{m}_{m-1}+\cdots+(-1)^m\mathfrak{m}_0=\beta_m-\beta_{m-1}+\cdots+(-1)^m\beta_0.
$$

An alternative approach to establish Morse theory for attractors of
nonsmooth systems  might be that one can still use Conley index
(which will not be developed here). The interested reader is
referred to Kunze,  Kupper and Li \cite{KKL} and Mrozek
\cite{Mrozek} etc. for Conley index theory of nonsmooth systems.

\section{Preliminaries}\label{prelim}
Let $X=\R^m$, which is equipped with the usual norm $|\.|$. For
convenience in statement, we will identify a single point $a\in X$
with the singleton $\{a\}$.

For any nonempty subsets $A$ and $B$ of $X$, define the {\em
Hausdorff semi-distance} and {\em Hausdorff distance}, respectively,
as

$$d_{\mbox{\tiny H}}(A,B)=\sup_{x \in A}d(x,B),\hs
\delta_{\mbox{\tiny H}}(A,B) = \max\left\{d_{\mbox{\tiny H}}(A,B),
d_{\mbox{\tiny H}}(B,A)\right\},
$$
where  $d(x,B)=\inf_{y\in B} |x-y|.$ We also assign
$d_{\tiny\mb{H}}(\emp,B)=0$.

  The {\em closure} of $A$ is denoted by $\ol A$, and the {\em interior} and {\em boundary} of $A$ are denoted by int$A$ and $\d A$, respectively. We
use $\cB(A,r)$ to denote the $r$-{\em neighborhood} of $A$,
i.e.,$$\cB(A,r)=\{y\in X|\,\,\,d(y,A)<r\}.$$ In particular,
$\cB_r=\cB(0,r)$ is the {\em ball} in $X$ centered at $0$ with
radius $r$.

We say that a subset $V$ of $X$ is a {\em neighborhood} of $A$, this
means $\ol A\subset \mbox{int}\,V$.

\subsection{Some basic  facts on differential inclusions}

Let $I$ be an interval.  A map $x(\cdot):\,I\rightarrow X$ is said
to be a {\em solution} of (\ref{e:1.1}) on $I$,  if it is absolutely
continuous on any compact interval $J\subset I$ and solves
(\ref{e:1.1}) at {\sl a.e.} $t\in I$.

A solution on $\R^1$ will be simply   called a {\em complete
solution}.

Let $x\in X$, and $A\subset X$. We denote by $\cS_x$ the family of
solutions $x(\.)$ of (\ref{e:1.1}) with initial value $x(0)=x$, and
$\cS_A=\cup_{x\in A} \cS_x$. Define set-valued map $\cR$ on $\R^+\X
X$ as:
 $$
\cR(t)x = \left\{x(t)|\,\,x(\cdot)\in \mathcal{S}_x \mb{ which
exists on }[0,t]\right\},\Hs \A \,(t,x)\in \R^+\X X.
$$Then
$\cR$ satisfies the following  semigroup property:
$$
\cR(0)x=x,\Hs \A \,x\in X,
$$
$$
\cR(s)\cR(t)x=\cR(s+t)x,\Hs \A\, s,t\geq0,\,\,x\in X.
$$
We will call $\cR$ the {\em multi-valued semiflow} (or {\em
set-valued semiflow}) generated by (\ref{e:1.1}), sometimes written
as $\cR(t)$.

\begin{lemma}\label{l:2.a1}  $\cite{AC,Deim}$.  Let $I$ be a compact interval. Then, for any sequence $\de_n$ $\ra$ $0$ and sequence
$x_n(\cdot)$ of absolutely continuous and uniformly bounded
functions on $I$ satisfying:
$$
x'_n(t)\in \ol{\mbox{con}}\
F\left(x_n(t)+\de_n\ol{\cB}_1\right)+\de_n\ol{\cB}_1,
$$
 there
exists a subsequence $x_{n_k}(\cdot)$ converging uniformly to a
solution $x(\cdot)$ of $(\ref{e:1.1})$ on $I$.
\end{lemma}

\begin{lemma}\label{l:2.a2} $\cite{Li2}$   Let $K$ be a compact subset of $X$
and let  $0<T<\8$. If  no solution $x(\cdot)\in \cS_K$  blows up on
$[0,T]$, then there exists an $R>0$ such that
$$
|x(t)|\leq R,\Hs \A t\in [0,T] ,\,\,x(\cdot)\in \cS_K.$$Consequently
$\cR([0,T])K$ is compact.
\end{lemma}

\subsection{Dynamical concepts and attractors in the strong sense}

As we emphasized in the introduction, all the dynamical concepts in
this work should be understood {\em in the strong sense}, except
otherwise statement.

 Let $A,B\subset X$. We say that
$A$ {\em attracts} $B$, this means that no solution
$x(\cdot)\in\cS_B$ blows up in finite time, moreover,
$$\lim_{t\ra+\8}d_{\mbox{{\tiny H}}}(\cR(t)B,\,A)=0.$$
The {\em attraction basin} $\W(A)$ of $A$ is defined  as:
$$\W(A)=\{x\in X|\,\,\mb{ $A$ attracts $x$}\}.$$

The set $A$ is said to be  {\em positively invariant} (resp. {\em
invariant}), if
$$ \cR(t)A\subset A\,\,\,\,(\,\mb{resp. }\cR(t)A=A\,),\Hs \A\, t\geq0.$$
$A$ is said to be {\em weakly invariant}, if for any $x\in A$, there
passes through $x$ a complete solution $x(\.)$ with $x(\R^1)\subset
A$. The {\em $\omega$-limit sets} $\omega(A)$ is defined as:
$$
\omega(A):=\{y\in X:\,\,\exists t_n\rightarrow\infty\mbox{ and }
y_n\in \cR(t_n)A \mb{ such that } y_n\rightarrow y\},
$$

\begin{definition}

Let $\a$ be a compact subset of $X$.  If there is a neighborhood $U$
of $\a$ such that $\a=\w(U)$, then we say that $\a$ is an {
attractor} of $(\ref{e:1.1})$ $($in the strong sense$)$.

We  allow the empty set $\emptyset$ to be an
 attractor with $\W(\emptyset)=\emptyset$.

  A { global attractor} is an
attractor $\a$ with $\W(\a)=X$.
\end{definition}

\bp\label{p:2.5}  $\cite{Li1,Li2}$ Let $\a$ be an attractor  of
$(\ref{e:1.1})$. Then
\begin{enumerate}\item[$(1)$] $\a$
is  invariant.
\item[$(2)$]
 $\a$ is  Lyapunov stable, that is, for any $\ve>0$, one can find a $\de>0$ such that
 $$
\cR(t)\mb{{\em\cB}}(\a,\de)\subset \mb{{\em\cB}}(\a,\ve),\Hs \A\,
t\geq 0.
$$
\item[$(3)$] $\W(\a)$ is a positively invariant open neighborhood of
$\a$.
\item[$(4)$] $\a$ attracts any compact subset $K$ of\, $\W(\a)$.\end{enumerate}

\ep

Concerning the existence of attractors, it is known that if there is
a compact set $K$ that attracts a neighborhood of itself, then
(\ref{e:1.1}) has an attractor $\a\subset K$ (see, e.g.,
\cite{Li2}).\vs Let there be given  an attractor $\a$ of
(\ref{e:1.1}).
\begin{definition}
An {\em attractor neighborhood} $\O$ of $\a$ means a positively
invariant neighborhood of $\a$ with $\O\subset
\W(\a)$.\end{definition}

Let $\a\subset U\subset X$. We define the {\em attraction basin}
$\W^U(\a)$ of $\a$ in $U$ as: \be\label{e:2.7}
\W^U(\a)=\W(\a)\cap\{x\in U|\,\,\cR(t)x\subset U\mb{ for all
}t\geq0\} . \ee
 One easily checks the validity of the following fact by using Lyapunov stability of $\a$.

\bp\label{p:2.2} Let $U$ be an open neighborhood of $\a$. Then
$\W^U(\a)$ is an open attractor neighborhood of $\a$. \ep

 Now we consider for $\de\geq 0$ the inflated system of
(\ref{e:1.1}):
\begin{equation}\label{e:1.2}
x'(t)\in \ol{\mb{con}}\ F\left(x(t)+\de\ol\cB_1\right)+\de\ol\cB_1
.\ee The proof of the robustness result below can be obtained by
slightly modifying the one for Theorem 2.10 in \cite{LK2} (see also
Theorem 2.9 in \cite{Lisiam}). We omit the details.

\begin{theorem}\label{t:2.10} Let $\a$ be an attractor of $(\ref{e:1.1})$. Then there exists a $\de_0>0$ such that the inflated
system $(\ref{e:1.2})$ has an attractor $\a(\de)$ for $\de<\de_0$.
Moreover, it holds that
\begin{enumerate}
\item[$(1)$] $\de_{\mbox{\em{\tiny H}}}\(\a(\de),\,\a\)\ra0$ as
$\de\ra0$;

\item[$(2)$] for any open neighborhood $U$ and compact set $K\subset
\W^U(\a)$, we have  $$K\subset\W^U(\a(\de)),$$ provided $\de$ is
 sufficiently small. \end{enumerate}
\end{theorem}

\subsection{Morse decompositions of attractors}

Let $\a$ be an attractor of (\ref{e:1.1}) with attraction basin
$\W=\W(\a)$.

 We say that a compact subset $A$ of $\a$ is an {\em attractor} in $\a$,
 this  means that there exists a neighborhood $U$ of $A$ such that $$\w\(U\cap
 \a\)=A.$$

\begin{theorem}\label{t:3.4}$\cite{Lisiam}$ Let $A\subset \a$ be an attractor in $\a$. Then $A$ is also an attractor
of $(\ref{e:1.1})$ in $X$.\end{theorem}

Let $A$ be an attractor in $\a$. Define $$A^*=\{x\in \a|\,\,
\w(x)\setminus A\neq\emptyset\}. $$
 $A^*$ is said to be a {\em repeller }of
(\ref{e:1.1}) in $\a$  dual to $A$, and $(A,A^*)$ is said to be an
{\em attractor-repeller pair} in $\a$.

It is known that $A^*$ is compact and weakly invariant
\cite{Lisiam}. Moreover, we have\be\label{e:2.3}A^*=\a\setminus
\W^\a(A)=\a\setminus \W(A).\ee

\begin{definition}\label{d:3.2}Let $\a$ be an attractor of
$(\ref{e:1.1})$. An ordered collection $\cM=\{M_1,\cdots,M_l\}$ of
compact subsets of $\a$ is called a Morse decomposition of $\a$, if
there exists an increasing sequence \be\label{e:2.4} \emp=A_0\subset
A_1\subset\cdots\subset A_l=\a \ee of attractors  in $\a$, called
{Morse filtration} of $\cM$, such that \be\label{e:2.6} M_k=A_k\cap
A_{k-1}^*,\Hs 1\leq k\leq l. \ee

The sets $M_k$ are called {Morse sets}.
\end{definition}
\br Since a Morse set $M_k$ may disappear under perturbations, for
convenience in statement, we allow $M_k$ to be the empty set $\emp$.
This occurs  in case $A_k=A_{k-1}$.
 However, if two Morse
decompositions $\cM$ and $\cM'$ have the same nonempty Morse sets,
they will be regarded as the same.\er

The following result taken from \cite{Lisiam}  provides some
information  on  Morse decompositions of attractors.

\begin{theorem}\label{t:3.2}

Let $\cM=\{M_1,\cdots,M_l\}$ be a Morse decomposition of $\a$ with
Morse filtration $ \emptyset=A_0\subset A_1\subset\cdots\subset
A_l=\a$. Then
\begin{enumerate}
\item[$(1)$] For each $k$, $(A_{k-1},M_k)$ is an
 attractor-repeller pair in $A_k$;

\item[$(2)$] $M_k$ are pair-wise disjoint weakly invariant compact sets;
 \item[$(3)$] If
$\gamma$ is a complete trajectory, then either $\gamma(\R)\subset
M_k$ for some Morse set $M_k$, or else there are indices $i<j$ such
that $ \alpha(\gamma)\subset M_j$ and $\w(\gamma)\subset M_i; $
\item[$(4)$] The attractors $A_k$ are uniquely determined by the
Morse sets, that is,
$$
A_k=\Cup_{1\leq i\leq k}W^u(M_i),\Hs 1\leq k\leq l,
$$
where $W^u(M_i)$ is the unstable manifold of $M_i$, namely, $$
W^u(M_i)=\{x|\,\,\mb{there is a trajectory $\gamma:\R\ra \a$ through
$x$ with $\alpha(\gamma)\subset M_i$}\}.$$

\end{enumerate}
\end{theorem}

\section{Smooth Lyapunov Functions of Attractors}

Let $U$ be an open subset of $X$.

  A nonnegative function
$\alpha\in C(U)$
  is said to be
{\em radially unbounded} on $U$, notated by $\alpha\in \cK^\8$,  if
for any $R>0$ there exists a compact subset $K\subset U$ such that
$$
\alpha(x)>R,\Hs \A\,x\in U\-K.
$$

A continuous function $V$ on an open neighborhood $U$ of an
attractor $\a$ is said to be a { \em Lyapunov function}  of $\a$ on
$U$, if \vs (1) $V$ is constant on $\a$;\vs (2) $V(x(t))$ is
strictly decreasing in $t$ for any solution $x(\.)$ of
$(\ref{e:1.1})$ in $U\-\a$. \vs

\noindent Our main results in this part are  the following two
theorems. One asserts the existence of  bounded smooth Lyapunov
functions of attractors defined on the whole phase space $X$, which
will be used to construct Morse-Lyapunov functions of attractors in
Section 4. The other relates to the existence of  smooth $\cK^\8$
Lyapunov functions of attractors, which  plays an important role in
discussing topological  properties of attractors.

\bt\label{t:l1}  Let $\O$ be an open attractor neighborhood of $\a$.
Then there exists a  smooth function $V\in C^\8(X)$ such that
\be\label{e:v15} V|_\a=0,\hs V|_{\a^c}>0,\hs V|_{\O^c}=1\ee
 \be\label{e:v16}\max_{v\in F(x)} \nab
V(x)\stackrel{\centerdot}{} v\leq -w(x),\Hs \A\,x\in X,\ee where
 $w\in C(X)$ is a nonnegative function
satisfying:\be\label{e:v17}w|_{\a\Cup \O^c}=0,\hs w|_{\O\-\a}>0. \ee
\et

\bt\label{t:l2} Let $\O$ be an open attractor neighborhood of $\a$.
Then there exists a $\cK^\8$ function $V\in C^\8(\O)$ that vanishes
on $\a$, such that \be\label{e:v2}\max_{v\in F(x)} \nab
V(x)\stackrel{\centerdot}{} v\leq -w(x),\Hs \A\,x\in\O,\ee where
 $w\in C(\O)$ is a nonnegative function
satisfying: \be\label{e:v1}w|_\a=0,\hs w|_{\O\-\a}>0.\ee \et

A smooth converse Lyapunov theorem for  strongly asymptotically
stable compact sets of (\ref{e:1.1}) in case $F(x)$ is continuous in
$x$ can be found in \cite{ST} (the approach therein is mainly based
on some arguments which seem to be more geometrical).

In what follows we first construct smooth Lyapunov functions when
$F$ is locally Lypschitz. Then we prove Theorem \ref{t:l1} by
appropriate approximations and some smoothing method used in
\cite{CL,LSW} etc. Theorem \ref{t:l2} follows directly from Theorem
\ref{t:l1}.

We remark  that although we follow the procedure in \cite{CL} and
\cite{LSW} etc., our construction method of the Lyapunov functions
here is quite different from those in the literature. In contrast,
ours seems to be more direct and simpler, and can be easily handled.

\subsection{Locally Lipschitz $F$ {vs} locally Lipschitz Lyapunov function}
In this subsection we assume that  $F$ is locally Lipschitz on $X$,
that is, for any compact subset $K\subset X$, there exists an $L>0$
such that \be F(x)\subset F(y)+L|x-y|\,\ol\cB_1,\Hs \A\,x,y\in K.
\ee

Let $U$ be an open subset of $X$. For $V\in C(U)$, $x\in U$ and
$v\in X$, define
$$
D^+_vV(x)=\limsup_{\tau\ra0^+}\frac{V(x+\tau v)-V(x)}{\tau}.$$
$D^+_vV(x)$ is called the {\em Dini supderivative} of $V$ along the
vector $v$.

The main result in this subsection is contained in the following
proposition.
 \begin{proposition}\label{p:2.1} Let $\O$ be an open attractor neighborhood of $\a$.
 Then there exists on $\O$ a locally Lipschitz  $\cK^\8$ function
 $V$ vanishing on $\a$ such that
\be\label{e:v2}\max_{v\in F(x)} D^+_vV(x)\leq -w(x),\Hs
\A\,x\in\O,\ee where
 $w\in C(\O)$ is a nonnegative function
satisfying $(\ref{e:v1})$.
\end{proposition}

The following lemma plays a basic role in the proof of Proposition
\ref{p:2.1}.

\bl\label{l:3.4} For any $\de>0$ with $\ol \cB(\a,\de)\subset \O$,
there exists  a  locally Lipschitz $\cK^\8$ function $V\in C(\O)$
such that \be\label{e:v13}V|_\a=0,\hs V|_{\O\-\ol{ \mb{\footnotesize
{\em B}}}(\a,\de)}>0,\ee \be\label{e:v14} \max_{v\in F(x)}
D^+_vV(x)\leq -V(x),\Hs x\in \O.\ee

\el

\noindent{\bf Proof.}  Pick a sequence of compact subsets $K_n$  of
$\O$ ($n=0,1,2\cdots$) such that $$ \ol\cB(\a,\de)\subset\mb{int}
K_0\subset K_0\subset \cdots\subset K_n\subset
\mb{int}K_{n+1}\cdots,\hs \O=\cup_{n\geq0}K_n.
$$
Then since the distance between $\d K_n$ and $K_{n-1}$ is positive,
there exists for each $n$ a $\tau_n>0$ such that
\be\label{e:3.1}\cR(t)\d K_n\subset \O\setminus K_{n-1},\Hs t\in
[0,\tau_n).\ee Take a nonnegative function $a_0\in C^\8(X)$ with
\be\label{e:2.2} a_0(x)\equiv 0\,\, \,(\mb{on }\cB(\a,\de/2)\,), \hs
\mb{and } a_0(x)>0\,\,\,(\mb{on }X\setminus \ol\cB(\a,\de)\,).\ee
For each $n\geq 2$ we choose a nonnegative function $a_n\in C^\8(X)$
with
$$
a_n(x)=\left\{\ba{ll}1, \hs\hs&x\in K_n\setminus K_{n-1},\\[1ex]
0, &x\in K_{n-2}\cup K_{n+1}^c.\ea\right.
$$
(Such functions can be easily obtained by appropriately smoothing
some continuous ones.\,) Let \be\label{e:2.1}
\alpha(x)=a_0(x)+\sum_{n=2}^{\8}\frac{n}{\tau_n}a_n(x). \ee Note
that the righthand side of (\ref{e:2.1}) is actually a finite sum
over $n$. Consequently $\alpha\in C^\8(X)$. We observe that
$$\alpha(x)\equiv0\,\,\,(\mb{on }\ol{\cB}(\a,\de/2)),\hs \alpha(x)>0\,\,\,(\mb{on }\ol{\cB}(\a,\,\de) ^c)\,.$$ Define $V(x)$ as:
$$
V(x)=\sup_{x(\.)\in \mathcal{S}_x}\int_0^\8e^t\alpha(x(t))\,dt.
$$
For each $K_n$, since $\a$ attracts $K_n$, there exists a $T>0$ such
that \be\label{e:3.2} \cR(t)K_n\subset \cB(\a,\de/2),\Hs t>T. \ee
This implies that $V$ is well defined.

One easily verifies by construction of $\alpha$ that
$$ V(x)\equiv 0\,\mb{ \,(on }\a),\Hs V(x)>0\,\mb{\, (on
}\O\-\ol\cB(\a,\de)\,).$$

We show that $V(x)$ is a $\cK^\8$ function on $\O$. Let  $x\in
\O\setminus K_n$ ($n\geq 2$). Take a solution $x(\.)\in \cS_x$, and
let
$$
s_n=\sup\{t>0|\,\,x([0,t))\subset \O\setminus K_{n-1}\}.
$$
Then $x(s_n)\in K_{n-1}$. By (\ref{e:3.1}) one necessarily has
$$
x(t)\in K_n\setminus K_{n-1},\Hs t\in [s_n-\tau_n,\,s_n),
$$
Therefore
$$  V(x)\geq \int_{s_n-\tau_n}^{s_n}\alpha(x(t))\,dt\geq \frac{n}{\tau_n}\int_{s_n-\tau_n}^{s_n}a_n(x(t))\,dt=n,
$$
and the conclusion follows.

We now prove that $V$ is locally Lipschitz.

Let $x\in K_n$. We claim that  there exists a $x_*(\.)\in\cS_x$ such
that \be\label{e:3.0} V(x)=\int_0^\8e^t\alpha(x_*(t))\,dt. \ee
Indeed,  by (\ref{e:3.2}) we find  that
$$
V(x)=\sup_{x(\.)\in\cS_x}\int_0^Te^t\alpha(x(t))\,dt.
$$
Let $x_k(\.)\in\cS_x$ be a sequence such that \be\label{e:3.15}
V(x)=\lim_{k\ra+\8}\int_0^Te^t\alpha(x_k(t))\,dt.\ee Then by
compactness there exists a subsequence of $x_k(\.)$, still denoted
by $x_k(\.)$, such that $x_k(\.)$ converges uniformly on $[0,T]$ to
a solution $x_*(\.)\in \cS_x $. Passing to the limit in
(\ref{e:3.15}) one immediately obtains the validity of
(\ref{e:3.0}).

Denote by $g(t,x)$ the unique closest point in $F(x)$ to $x_*'(t)$.
Then $g:\R^1\X X\ra X$ is a Carath$\acute{\mb{e}}$odory function
(see \cite{Deim}, pp. 49), and therefore solutions of the ODE:
\be\label{e:3.16} y'(t)=g(t,y(t)) \ee exist at least locally. For
each $y\in K_n$, we pick  a solution $y(\.)$ of (\ref{e:3.16}).
Clearly $y(\.)\in \cS_y\, $, and hence it exists for all $t\geq0$.
As $\cR([0,T])K_n$ is compact, by virtue of (\ref{e:3.2}) we deduce
that there exists an $N>n$ such that
$$
\cR(t)K_n\subset K_N,\Hs \A\,t\geq0.
$$
Since $F$ is locally Lipschitz on $X$, there is an $L>0$ such that
$$
F(x_1)\subset F(x_2)+L|x_1-x_2|\ol\cB_1,\Hs \A\,x_1,x_2\in K_N.
$$
Thus
$$\ba{ll}|y'(t)-x_*'(t)|&=|x_*'(t)-g(t,y(t))|=d\(x_*'(t),\,F(y(t))\)\\[1ex]&\leq
(\mb{recall }x_*'(t)\in F(x_*(t))\,)\leq L|x_*(t)-y(t)|,\Hs t\geq
0,\ea$$ which implies
$$
|y(t)-x_*(t)|\leq |x-y|e^{Lt},\Hs t\geq 0.
$$
Therefore
$$\ba{lll}
V(x)&=\int_0^\8e^t\alpha(x_*(t))\,dt=\int_0^Te^t\alpha(x_*(t))\,dt\\[2ex]
&=\int_0^Te^t\alpha(y(t))\,dt+\int_0^Te^t\(\,\alpha(x_*(t))-\alpha(y(t))\,\)\,dt\\[2ex]
&\leq V(y)+M|x-y|\int_0^Te^{(L+1)t}dt,
 \ea
$$
where $M>0$ is the Lipschitz constant of $\alpha$ on $K_N$. Because
$x,y\in K_n$ are arbitrary, we conclude that
$$
|V(x)-V(y)|\leq M\int_0^Te^{(L+1)t}dt\,|x-y|,\Hs \A\, x,y\in K_n.
$$

Now let us evaluate $D_v^+V(x)$ for $x\in \O$ and $v\in F(x)$.

 For every $y$ denote by $g(y)\in F( y)$ the unique closest point in $F( y)$ to $v$. Then
 the function $g$ is continuous (see \cite{AC}), and  $ g(x)=v$. Let $x(\.)$ be a solution to the
following initial value problem:
$$
x'(t)=g(x(t)),\Hs x(0)=x.
$$
Of course $x(\.)\in \cS_x $. Note that
\be\label{e:3.17}x(\tau)=x+\tau v+o(\tau).\ee We
write\be\label{e:3.5} \frac{V(x+\tau
v)-V(x)}{\tau}=\frac{V(x(\tau))-V(x)}{\tau}+\frac{V(x+\tau
v)-V(x(\tau))}{\tau}. \ee Since $V$ is locally Lipschitz, by
(\ref{e:3.17}) one easily checks that the second term in the
righthand side of the above equation goes to $0$ as $\tau\ra0$. To
estimate  the first term in the righthand side, we first observe
that
$$
\ba{lll}V(x(\tau))&=\sup_{y(\.)\in\cS_{x(\tau)}}\int_0^{+\8}e^t\alpha(y(t))\,dt\\[2ex]
&=\sup_{y(\.)\in\cS_{x(\tau)}}e^{-\tau}\int_0^{+\8}e^{t+\tau}\alpha(y(t))\,dt.\ea
$$
For any $y=y(\.)\in \cS_{x(\tau)}$, define
$$
x_y(t)=\left\{\ba{ll}x(t),\Hs & t\in[0,\tau];\\[1ex]
y(t-\tau),\hs\hs & t>\tau.\ea\right.
$$
Then $x_y(\.)\in\cS_x $. We can rewrite $V(x(\tau))$ as
$$\ba{ll}
V(x(\tau))&=\sup_{y(\.)\in\cS_{x(\tau)}}e^{-\tau}\int_0^{+\8}e^{t+\tau}\alpha(x_y(t+\tau))\,dt\\[2ex]
&=\sup_{y(\.)\in\cS_{x(\tau)}}e^{-\tau}\int_\tau^{+\8}e^{t}\,\alpha(x_y(t))\,dt,\ea$$
from which it can be easily seen that
$$\ba{ll}
V(x(\tau))&\leq\sup_{x(\.)\in\cS_x
}e^{-\tau}\int_\tau^{+\8}e^{t}\,\alpha(x(t))\,dt\leq
e^{-\tau}V(x),\ea$$ which implies
$$
\limsup_{\tau\ra0^+}\frac{V(x(\tau))-V(x)}{\tau}\leq -V(x).
$$
By (\ref{e:3.5}) one  immediately deduces that
$$
D^+_vV(x)\leq -V(x).
$$

The proof of the lemma is complete. $\Box$

\Vs

\noindent{\bf Proof of Proposition \ref{p:2.1}}. Take a decreasing
sequence of numbers $\de_n\downarrow0$ with
$$\ol\cB(\a,2\de_0)\subset\O.$$ Then for each $\de_n$  one can find a
$\cK^\8$  function $V_n\in C(\O)$ satisfying all the properties in
Lemma \ref{l:3.4} with $V$ and $\de$ therein replaced by $V_n$ and
$\de_n$, respectively. In particular, $V_n$ vanishes on $\a$ and
satisfies:
$$
V_n(x)>0,\Hs\mb{if } x\not\in \ol\cB(\a,\de_n).
$$

For each $n$ we  pick an $R_n>0$ so that $$ \ol\cB(\a,2\de_0)\subset
K_n:=\{x\in\O|\,\,V_n(x)\leq R_n\}.$$ Define \be\label{e:2.5}
W_n(x)=\left\{\ba{ll} V_n(x),\hs \hs&x\in K_n;\\[1ex]
R_n,&x\in X\setminus K_n.\ea\right. \ee Then $W_n$ is globally
Lipschitz on $X$. Denote by $L_n$ the Lipschitz constant of $W_n$.
Let
$$
V(x)=V_0(x)+\sum_{n=1}^\8\frac{1}{2^nR_nL_n}W_n(x),\Hs x\in\O.
$$
One trivially checks that $V$ is a locally Lipschitz $\cK^\8$
function on $\O$. Moreover, $V|_\a=0$.

Now we evaluate $D^+_vV(x)$ for any $x\in \O$ and $v\in F(x)$.
Observe (by (\ref{e:v14})) that \be\label{e:3.20}
D^+_vW_n(x)\leq\left\{\ba{ll} -V_n(x),\hs\hs &\mb{if }x\in\mb{int}K_n;\\[2ex]
0,&\mb{otherwise }.\ea\right. \ee So if we take $w_n\in C(\O)$ as:
$$
w_n(x)=\min\(V_n(x),\,d(x,K_n^c)\),\Hs x\in \O,
$$
then by (\ref{e:3.20}) one finds that
$$
D^+_vW_n(x)\leq -w_n(x). $$ Note that $0\leq w_n(x)\leq R_n$ for all
$x\in\O$. Define
$$
w(x)=V_0(x)+\sum_{n=1}^\8\frac{1}{2^nR_nL_n}w_n(x),\Hs x\in\O.
$$
Then $w\in C(\O)$. We verify that
$$
w(x)>0,\Hs \A\,x\in\O\-\a.
$$
There are two possibilities.

(1)  ``$x\not\in \ol\cB(\a,\de_0)$\,''. In this case we directly
have
$$w(x)\geq V_0(x)>0.$$

(2) ``$x\in\ol\cB(\a,\de_0)$\,''.  We pick  an $n$ large enough so
that $x\not\in\ol\cB(\a,\de_n)$. Then
$$x\in \ol\cB(\a,\de_0)\-\ol\cB(\a,\de_n)\subset
\mb{int}K_n\-\ol\cB(\a,\de_n).$$ Since both $V_n(x)$ and $d(x,
K_n^c)$ are positive, by the definition of $w_n$ we deduce that
$w_n(x)>0$. Therefore
$$
w(x)\geq\frac{1}{2^nR_nL_n}w_n(x)>0.
$$

Now it is easy to deduce  that
$$
D^+_vV(x)\leq
D^+_vV_0(x)+\sum_{n=1}^\8\frac{1}{2^nR_nL_n}D^+_vW_n(x)\leq -w(x),
$$
which  also follows
$$
V(x)>0,\Hs \A\,x\in\O\-\a.
$$

The proof is complete. $\Box$

\subsection{Locally Lipschitz $F$ vs smooth Lyapunov function}
 \begin{proposition}\label{p:3.5} Let $\O$ be an open attractor neighborhood of $\a$. Then there exists a
$\mathcal{K}^\8$ function $V\in C^\8(\O)$ which vanishes on $\a$,
such that \be\label{e:v18}\max_{v\in F(x)}\nab
V(x)\stackrel{\centerdot}{} v\leq -w(x),\Hs \A\,x\in\O,\ee where
 $w\in C(\O)$ is a nonnegative function satisfying $(\ref{e:v1})$.
\end{proposition}

\Vs\noindent{\bf Proof.} Following the procedure  as in \cite{CL},
Section 5 (with very minor modifications), one can obtain a Lyapunov
function $V$ of $\a$, which belongs to  $C^\8(\O\-\a)$ and satisfies
(\ref{e:v18}),  by smoothing the Lipschitz one given in Proposition
\ref{p:2.1} above.
 Further applying the following
lemma, which is a slightly modified version of Lemma 4.3 in
\cite{LSW}, we immediately get a smooth Lyapunov function on $\O$,
as desired. $\Box$

\bl Assume that $V\in C(\O)\Cap C^\8(\O\-\a),$ and that
$$
V|_\a=0,\hs V|_{\O\-\a}>0.
$$
Then there exists a nonnegative function $\beta\in C^\8([0,\8))$
with
$$
\beta^{(k)}(0)=0\,\,\,(k=0,1,\cdots),\hs \lim_{t\ra+\8}\b(t)=+\8
$$
and
$$
\b'(t)>0,\Hs \A\,t>0,
$$such that $W(x)=\b(V(x))$ is a $C^\8$ function on $\O$.

\el

\noindent{\bf Proof.} See Lemma 4.3 in \cite{LSW}.  $\Box$

\subsection{The general case: proofs of Theorems \ref{t:l1} and \ref{t:l2}}

To prove Theorems \ref{t:l1} and \ref{t:l2}, we need the following
approximation lemma. Closely related results can also be found  in
Lasry-Robert \cite{Lasry} etc.

  \bl\label{l:3.8} For any $\de,R>0$, one can find a Lipschitz continuous  multifunctions $F_L: X\mapsto
\mathcal{C}(X)$ such that \be\label{e:3.9}F(x)\subset F_L(x)\subset
\ol{\mb{\em con}}\,F(x+\de\ol B_1)+\de\ol\cB_1,\Hs \A\,x\in
\ol\cB_R,\ee where $\mathcal{C}(X)$ denotes the family of nonempty
compact convex subsets of $X$.

 \el

\noindent{\bf Proof.} The proof  is directly adapted from that for
Theorem 9.2.1 in \cite{AF}.

Let $\de,R>0$ be given arbitrary. Then for every $x\in X$, by upper
semicontinuity of $F$ there exists $0<r_x<\de$ such that
$$
F(\ol\cB(x,r_x))\subset F(x)+\de\ol\cB_1.
$$
The family of the balls $\cB(x,r_x/4)$ ($x\in \ol\cB_R)$  cover
$\ol\cB_R$. So there exist a finite number of balls
$\cB(x_i,r_{x_i}/4)$ ($1\leq i\leq n$)  that cover $\ol\cB_R$.

We rewrite  $r_{x_i}$ as $r_i$ for simplicity and take a Lipschitz
partition of unity $\{a_i(x)\}_{1\leq i\leq n}$ subordinated to this
finite covering. That is, each $a_i:X\ra[0,1]$ is a Lipschitz
function that vanishes outside $\cB(x_i,r_i/4)$, and
$$
\sum_{1\leq i\leq n}a_i(x)=1,\Hs x\in \ol\cB_R.
$$
Define $F_L$ as:
$$
F_L(x)=\sum_{1\leq i\leq n}a_i(x)\,F\(\ol\cB(x_i,r_i/4)\),\Hs x\in
X.
$$
Then $F_L$ is Lipschitz.

Let $x\in\ol\cB_R$, and let $I(x)$ be the set of all $i$ such that
$a_i(x)\ne0$. Then we must have $$x\in \cB(x_i,r_i/4),\Hs \A\, i\in
I(x).$$ Therefore for all $i,j\in I(x)$, \be\label{e:3.10}
|x_i-x_j|\leq|x_i-x|+|x-x_j|\leq r_i/2. \ee Fix a $k\in I(x)$
satisfying
$$
r_k=\max_{i\in I(x)}r_i.
$$
Then $|x-x_k|\leq r_k/4$, and (\ref{e:3.10}) implies
$$
x_i\in \cB(x_k,r_k/2),\Hs\A\, i\in I(x).
$$
Consequently for all $i\in I(x)$,
$$\ba{ll}
F(\ol\cB(x_i,r_i/4))&\subset F(\ol\cB(x_k,r_k))\subset
F(x_k)+\de\ol\cB_1\\[1ex]
&\subset F(\ol\cB(x,r_k))+\de\ol\cB_1\subset
F(\ol\cB(x,\de))+\de\ol\cB_1\subset \ol{\mb{con}}\,
F(x+\de\ol\cB_1))+\de\ol\cB_1 .\ea
$$
Since the last term of this relation is convex, from the definition
of $F_L$ we see that
$$
F_L(x)\subset \ol{\mb{con}}\, F(x+\de\ol\cB_1))+\de\ol\cB_1.
$$

There remains to verify that $F(x)\subset F_L(x)$. Assume that $y\in
F(x)$. Then
$$
y\in F(x)\subset F(\ol\cB(x_i,r_i/4)),\Hs \A\, i\in I(x).
$$
Therefore
$$
y=\sum_{i\in I(x)}a_i(x)y\in \sum_{i\in
I(x)}a_i(x)F(\ol\cB(x_i,r_i/4))=F_L(x).
$$

The proof is complete. $\Box$

\Vs Now we are in a position to prove  Theorems \ref{t:l2} and
\ref{t:l1}.\Vs

\noindent{\bf Proof of Theorem  \ref{t:l1}.} Let $\mathcal{O}$ be an
open attractor neighborhood
  of $\a$. Take a sequence of
compact neighborhoods $K_n$ of $\a$ with
$$
K_1\subset\mb{int}K_2\subset \cdots\subset K_n\subset
\mb{int}K_{n+1}\cdots, \hs \mathcal{O}=\cup_{n\geq1}K_n.$$ For each
$K_n$, since $\a$ attracts $K_n$ we deduce that $\cR(t)K_n$ is
bounded for $t\geq0$. Hence there exists an $R_n>0$ such that
$$
\cR(t)K_n\subset \cB_{R_n},\Hs t\geq 0.
$$
Let $U_n=\O\cap \cB_{2R_n}$. Then clearly
$$
K_n\subset \W^{U_n}(\a)\subset \O,
$$where $\W^{U_n}(\a)$ is the attraction basin of $\a$ in $U_n$; see
(\ref{e:2.7}) for the definition. By virtue of Theorem \ref{t:2.10}
there exists $\de>0$ such that the inflated system (\ref{e:1.2}) has
an attractor $\a(\de)$ with
$$
\a(\de)\subset\mb{int}\, K_n\cap \cB(\a,1/n),\hs
K_n\subset\W^{U_n}(\a(\de)).$$

Thanks to Lemma \ref{l:3.8} we can take a Lipschitz function $F_L:
X\mapsto\cC(X)$ such that \be\label{e:3.11} F(x)\subset
F_L(x)\subset \ol{\mb{con}}\,F(x+\de\ol\cB_1)+\de\ol\cB_1,\Hs
\A\,x\in \ol\cB_{3R_n}.\ee Consider the system \be\label{e:3.25}
x'(t)\in F_L(x(t)).\ee The second inclusion in (\ref{e:3.11})
implies that (\ref{e:3.25}) has an attractor, denoted by  $\a_n$. We
infer from (\ref{e:3.11})  that
$$
\a\subset\a_n\subset \a(\de). $$ We claim that \be\label{e:3.26}
\W^{U_n}(\a_n)\subset \W^{U_n}(\a)\subset\O.\ee Indeed, since
$K_n\subset \W^{U_n}(\a(\de))\cap \W^{U_n}(\a)$, it is easy to see
that
$$K_n\subset \W^{U_n}(\a_n).$$ Assume that $x\in
\W^{U_n}(\a_n)$. Let $\cR_{F_L}(t)$ be the multi-valued semiflow
generated by (\ref{e:3.25}). Then by
$\a_n\subset\a(\de)\subset\mb{int}K_n$ we deduce that
$\cR_{F_L}(t)x\subset \mb{int}K_n$ for $t$ sufficiently large. It
then follows that \be\label{e:3.27}\cR(t)x\subset
\cR_{F_L}(t)x\subset \mb{int}K_n\subset K_n\ee for $t$ sufficiently
large. Because $\a$ attracts $K_n$, we conclude by (\ref{e:3.27})
that $\a$ attracts $x$. Note that $x\in \W^{U_n}(\a_n)$ also implies
$$
\cR(t)x\subset \cR_{F_L}(t)x\subset U_n, \Hs\mb{for all } t\geq 0.$$
Hence $x\in  \W^{U_n}(\a)$. This proves (\ref{e:3.26}).

$\O_n:=\W^{U_n}(\a_n)$ is an open attractor neighborhood of $\a_n$.
Therefore Proposition \ref{p:3.5} implies that there exists  a
$\cK^\8$ function $V_n\in C^\8(\O_n)$ which vanishes on $\a_n$ such
that \be\label{e:v20}\max_{v\in F_L(x)}\nab
V_n(x)\stackrel{\centerdot}{} v\leq -w_n(x),\Hs \A\,x\in\O_n,\ee
where
 $w_n\in C(\O_n)$ is a nonnegative function satisfying:
 \be\label{e:v21}w_n|_{\a_n}=0,\hs w_n|_{\O_n\-{\a_n}}>0.\ee

For each $n$ take an $a_n>0$ sufficiently large so that
$$
K_n\subset D_n:=\{x\in\O_n|\,\, V_n(x)\leq a_n\}\subset\O.
$$ Define
$$
\psi_n(x)=\left\{\ba{ll} V_n(x),\hs \hs&x\in D_n;\\[1ex]
a_n,&x\in X\setminus D_n.\ea\right. $$ Clearly $\psi_n\in
C^\8(X\-\pa D_n)$. We will make a slight modification to $\psi_n$ to
obtain a smooth function $ \phi_n\in C^\8(X)$. For this purpose we
first note that $\psi_n$ is globally Lipschitz on $X$. Therefore
there exists a $C>0$ such that \be\label{s1}
|\psi_n(y)-\psi_n(x)|\leq C|x-y|,\Hs \A\,x\in\pa D_n,\,\, y\in X.\ee
Set $$ G(s)=\left\{\ba{lll}\mb{sgn}\,(s) e^{-1/s^2},\hs\hs& s>0;\\[1ex]
0,& s=0,\ea\right.$$ where sgn($\.$) is the signal function. Then
$G\in C^\8(\R^1)$, and
$$
G^{(n)}(0)=0\,\,\,(\A\, n\geq0),\Hs G'(s)>0\,\,(\,\A\,x\ne0\,).
$$
Define
$$
\phi_n(x)=G\(\psi_n(x)-a_n\)+G(a_n).
$$
By (\ref{s1}) one easily checks that $\phi_n\in C^\8(X)$ with $$
\frac{\pa^l\phi_n(x)}{\pa x_{i_1}\cdots\pa x_{i_l}}=0,\Hs\A\,x\in\pa
D_n,\,\,\,\,l=1,2,\cdots,$$where $1\leq i_k\leq m$. Clearly $\phi_n$
satisfies: \be\label{e:3.29} \phi_n|_{\a_n}=0,\hs
\phi_n|_{\a_n^c}>0,\hs \phi_n|_{D_n^c}= G(a_n).\ee

For $x\in \mb{int}D_n$ and $v\in F(x)\subset F_L$, we have
\be\label{e:3.28}\ba{ll}\nab\phi_n(x)\stackrel{\centerdot}{} v&=G'\(\psi_n(x)-a_n\)\nab \psi_n(x)\stackrel{\centerdot}{} v\\[1ex]
&=G'\(\psi_n(x)-a_n\)\nab V_n(x)\stackrel{\centerdot}{} v \leq -
G'\(\psi_n(x)-a_n\)w_n(x).\ea \ee Let $$
\widetilde{w}_n(x)=\left\{\ba{ll}G'\(\psi_n(x)-a_n\)w_n(x),\hs\hs
&x\in
D_n\,;\\[1ex]
0,&x\in D_n^c\,.\ea\right.
$$
Since $D_n$ is a compact subset of $\O_n$ and $w_n\in C(\O_n)$, one
easily sees that  $\~w_n$ is  continuous and bounded on $X$.
(\ref{e:3.28}) amounts to say that
\be\label{e:3.14}\nab\phi_n(x)\stackrel{\centerdot}{} v\leq -
\widetilde{w}_n(x),\Hs \A\,x\in \mb{int}D_n,\,\,\,v\in F(x). \ee The
above estimate naturally holds  if $x\in X\-\mb{int}D_n$, as in this
case both sides vanish.

We also note that \be\label{e:3.13} \widetilde{w}_n(x)>0,\Hs \A\,
x\in\mb{int}D_n\-\a_n. \ee

$\phi_n$ is constant on $D_n^c$\,, so  $
||\phi_n||_{C^k(X)}=||\phi_n||_{C^k(D_n)} $ for all integers
$k\geq0$. Let
$$c_n=||\phi_n||_{C^n(X)}+||\widetilde{w}_n||_{C(X)}+1.$$ Define
\be\label{e:3.12}
V(x)=\gamma\sum_{n=1}^\8\frac{1}{2^nc_n}\phi_n(x),\Hs x\in X, \ee
where $\gamma=1/\sum_{n=1}^\8\frac{1}{2^nc_n}G(a_n)$. For any $ l\in
\mathbb{N}$ and $1\leq i_1,\cdots,i_l\leq m$, because
$$\sum_{n=l}^\8\frac{1}{2^nc_n}\left|\frac{\pa^l\phi_n(x)}{\pa
x_{i_1}\cdots\pa x_{i_l}}\right|  \leq
\sum_{n=l}^\8\frac{1}{2^nc_n}||\phi_n||_{C^l(X)}\leq
\sum_{n=l}^\8\frac{1}{2^nc_n}||\phi_n||_{C^n(X)}\leq
\sum_{n=l}^\8\frac{1}{2^n},
$$ we deduce that the series $\sum_{n=1}^\8\frac{1}{2^nc_n}\,\frac{\pa^l\phi_n(x)}{\pa
x_{i_1}\cdots\pa x_{i_l}}$ is uniformly convergent on $X$. It then
follows that
  $V\in C^\8(X)$. Moreover,
for any $x\in X$ and $v\in F(x)$, we have \be\label{e:3.30} \nab
V(x)\stackrel{\centerdot}{} v\leq
\gamma\sum_{n=1}^\8\frac{1}{2^nc_n}\nab\phi_n(x)\stackrel{\centerdot}{}
v\leq -w(x), \ee where $$w(x)=\gamma\sum_{n=1}^\8\frac{1}{2^nc_n}
\widetilde{w}_n(x).$$

There remains to check that $V$ and $w$ satisfies all the other
properties required in the Theorem.

Let $x\in \O^c$. Then $x\in D_n^c$ for all $n\in \mathbb{N}$. By the
construction of $V$ and $w$, we have
$$
V(x)=\gamma\sum_{n=1}^\8\frac{1}{2^nc_n}\phi_n(x)=\gamma\sum_{n=1}^\8\frac{1}{2^nc_n}G(a_n)=1,
$$
$$w(x)=\gamma\sum_{n=1}^\8\frac{1}{2^nc_n}
\widetilde{w}_n(x)=0.$$

Assume that $x\in \O\-\a$. In this case one can find a $j\in
\mathbb{N}$ sufficiently large so that
$$
x\in K_{j-1}\-\a_j\subset \mb{int}K_j\-\a_j\subset
\mb{int}D_j\-\a_j.
$$
(Recall that $\mathcal{O}=\cup_{n\geq1}K_n$ and $\a_n\subset
\cB(\a,1/n)$.) (\ref{e:3.13}) then implies that
\be\label{e:3.31}w(x)=\gamma\sum_{n=1}^\8\frac{1}{2^nc_n}
\widetilde{w}_n(x)\geq \gamma\frac{1}{2^jc_j} \widetilde{w}_j(x)
>0.\ee

Finally it is trivial to examine that both $V$ and $w$ vanish on
$\a$. That $V$ is positive on $\W\-\a$ directly follows from
(\ref{e:3.30}) and (\ref{e:3.31}).

The proof of the lemma is complete.
 $\Box$

\Vs \noindent{\bf Proof of Theorem \ref{t:l2}.} Let $V$ be a
Lyapunov function of $\a$ given by Theorem \ref{t:l1}. Define
$$L(x)=\eta(V(x)),\Hs x\in\O,$$ where
 $\eta(s)=-\ln(1-s)$ ($s\in[0,1)$). Then $L$ is a $\cK^\8$
Lyapunov function of $\a$ on $\O$ that satisfies all the desired
properties in the theorem. $\Box$

\section{Smooth Morse-Lyapunov Functions}

Let  $\a$ be an attractor of (\ref{e:1.1}) with  Morse decomposition
$\cM=\{M_1,\cdots,M_l\}$, and let $$ \cD=\cup_{1\leq k\leq l}M_k.
$$

 A continuous function $V$ on the attraction basin $\W=\W(\a)$
is said to be a {\em Morse-Lyapunov function} $($M-L function in
short\,$)$ of $\cM$ on $\W$, if

\Vs$(1)$ $V$ is constant on each Morse set $M_k$; \vs $(2)$
$V(x(t))$ is strictly decreasing in $t$ for any solution $x(\.)$ of
$(\ref{e:1.1})$ in $\W\-\cD$.\Vs

An M-L function $V$ of $\cM$ is said to be a {\em strict M-L
function}, if in addition it satisfies:
$$V(M_i)<V(M_j)$$
whenever $i<j$ with $M_i\ne\emp\ne M_j$\,.

The main result in this section is contained in the following
theorem.

\bt\label{t:3.1} $(${\bf Existence of smooth M-L functions}$)$ $\cM$
has a radially  unbounded strict M-L function $V\in C^\8(\W)$ such
that \be\label{e:v8}\max_{v\in F(x)} \nab
V(x)\stackrel{\centerdot}{} v\leq -w(x),\Hs\A\, x\in \W,\ee where
 $w\in C(\W)$ is a nonnegative function
satisfying \be\label{e:3.7} w|_\cD= 0,\hs w|_{\W\-\cD}>0. \ee

\et

\noindent{\bf Proof.} Let $$ \emp=A_0\subset A_1\subset\cdots\subset
A_l=\a$$ be the Morse filtration of $\cM$.

For $k=l$, by Theorem \ref{t:l2}  there exists a $\cK^\8$ function
$V_l\in C^\8(\W)$ such that $$V_l|_{A_l}=V_l|_{\a}=0,$$
$$\max_{v\in F(x)} \nab V_l(x)\cdt v\leq -w_l(x),\Hs\A\, x\in \W,$$
where
 $w_l\in C(\W)$ is a nonnegative function
satisfying $$ w_l|_{\a_l}=0,\hs w_l|_{\W\-\a_l}>0.$$

For each $1\leq k\leq l-1$,  it follows from Theorem \ref{t:l1} that
there exists a $V_k\in C^\8(X)$ such that $$V_k|_{A_k}=0, \hs
V_k|_{\W(A_k)^c}=1,$$ and $$\max_{v\in F(x)} \nab V_k(x)\cdt v\leq
-w_k(x),\Hs \A\,x\in X,$$ where
 $w_k\in C(X)$ is a nonnegative function
with $$w_k|_{A_k\cup\W(A_k)^c}=0,\hs w_k|_{\W(A_k)\setminus A_k}>0.
$$ Define $V\in C(\W)$ as: $$ V(x)=\sum_{1\leq k\leq
l}V_k(x), \Hs x\in\W. $$ We show that $V$ has all the required
properties with
$$w(x)=\sum_{1\leq k\leq l}w_k(x).$$

Let $x\in\W$, and $v\in F(x)$. Then
$$
\nab V(x)\stackrel{\stackrel{\centerdot}{}}{} v\leq \sum_{1\leq
k\leq l}\nab V_k(x)\stackrel{\centerdot}{} v\leq -\sum_{1\leq k\leq
l}w_k(x)=-w(x).
$$
This verifies (\ref{e:v8}). In what follows we examine the validity
of (\ref{e:3.7}).

Let $x\in\cD$. We may assume that $x\in M_k=A_k\cap A_{k-1}^*$ for
some $k$. If $i\geq k$, then we have $x\in A_i$, and hence
$w_i(x)=0$. On the other hand, in case  $i\leq k-1$ we deduce by
$x\in A_{k-1}^*$ that
$$x\in \W(A_{k-1})^c\subset \W(A_i)^c,$$ which also  implies $w_i(x)=0$. In
conclusion, $$w_i(x)=0,\Hs \mb{for all }1\leq i\leq l.$$
Consequently $w(x)=0$.

Now assume $x\in \W\setminus \mathcal{D}$. Then there is a smallest
$k$ such that $x\in \W(A_k)$. We claim that $x\not\in A_k$. Indeed,
if $x\in A_{k}$, then since  $(A_{k-1},M_k)$ is an
 attractor-repeller pair in $A_k$, we have
either $x\in M_k$, or $x\in\W(A_{k-1})$. The former case contradicts
to that $x\in \W\setminus \mathcal{D}$, and the latter one
 to the definition of $k$. Hence the claim holds true.
Now that $x\in \W(A_k)\-A_k$, one finds that  $w_k(x)>0$. Thus we
have
$$w(x)\geq w_k(x)>0.$$

Finally we check that $V$ is a strict M-L function of $\a$. Recall
that
$$
M_k=A_k\cap
A_{k-1}^*=A_k\cap\(\a\setminus\W(A_{k-1})\)=A_k\cap\W(A_{k-1})^c.
$$
 Thus if $i\leq k-1$,
then $ M_k\subset  \W(A_{k-1})^c\subset \W(A_i)^c, $ which implies,
in case $M_k\ne\emp$, that
$$V_i(M_k)=1,\Hs\mb{for } 1\leq i\leq k-1.$$  On the other hand if $i\geq k$, then
$M_k\subset A_k\subset A_i$. Therefore
$$V_i(M_k)=0,\Hs\mb{for }  i\geq k.$$ Hence if $M_k$ is nonvoid, then we conclude that
$$
V(M_k)=\sum_{1\leq i\leq l}V_i(M_k)=k-1.
$$

 The proof is complete.  $\Box$

\section{Morse Theory}

In this section we try to establish   a Morse theory for attractors.
First, we  prove  that all the  open attractor neighborhoods of an
attractor have the same homotopy type by employing smooth Lyapunov
functions. Then, based on this fact we introduce the concept of
critical groups for Morse sets and establish Morse equations and
inequalities.

For any $a\in\R^1$, we will denote by $V_a$ the level set of a
function $V:R^m\ra\R^1$,
$$
V_a=\{x\in\R^m|\,\,V(x)\leq a\}.
$$

\subsection{Homotopy equivalence of open attractor neighborhoods}
Let $\a$ be an attractor of the system (\ref{e:1.1}). We shall prove
that all the open  neighborhoods of $\a$ have the same homotopy
type.

The following lemma will play an important role. \bl\label{l:6.1}
Let $\O_i$ $(i=1,2)$ be two open attractor neighborhoods of $\a$,
and let $V_i$ $(i=1,2)$ be  Lyapunov functions of  $\a$ on $\O_i$
given by Theorem \ref{t:l2}, respectively.

Then for any compact set $K\subset \(\O_1\cap\O_2\)\-\a$, there
exists a smooth vector field $\Psi$ defined on $X$ such that $$ \nab
V_i(x)\stackrel{\centerdot}{} \Psi(x)<0,\Hs \A\, x\in K,\,\,\,i=1,2.
$$ \el

\noindent {\bf Proof.} Clearly  $\O:=\(\O_1\cap\O_2\)\-\a$ is open.
For each $x\in \O$ we fix a $v_x\in F(x)$. Then
$$
\nab V_i(x)\stackrel{\centerdot}{} v_x<0,\Hs i=1,2.
$$
Take an $r_x>0$ sufficiently small so that $\cB(x,r_x)\subset \O$,
and
$$\nab V_i(y)\stackrel{\centerdot}{} v_x<0,\Hs \A\,y\in \cB(x,r_x),\,\,\,i=1,2.
$$
Then the family of balls $\{\cB(x,r_x)\}_{x\in K}$ forms an open
covering of $K$, therefore by compactness of $K$ it has a finite
subcovering $\cV=\{\cB(x_k,r_{x_k})\}_{1\leq k\leq n}$\,. Let
$a_k\in C^\8(X)$ ($1\leq k\leq n$) be a smooth unit partition of $K$
subordinated to $\cV$, namely, each $a_k$ vanishes outside
$\cB(x_k,r_{x_k})$, and
$$
\sum_{1\leq k\leq n}a_k(x)\equiv 1,\Hs x\in K.
$$
Define $\Psi$ on $X$ as:
$$
\Psi(x)=\sum_{1\leq k\leq n}a_k(x)v_{x_k}\,,\Hs x\in X.
$$
Then $\Psi\in C^\8$. For each $x\in K$, we have
$$\ba{ll}
\nab V_i(x)\stackrel{\centerdot}{} \Psi(x)&=\sum_{1\leq k\leq
n}\,a_k(x)\,\(\nab
V_i(x)\stackrel{\centerdot}{} v_{x_k}\)\\[2ex]&\leq \max_{1\leq k\leq n}\(\nab
V_i(x)\stackrel{\centerdot}{} v_{x_k}\)\sum_{1\leq k\leq n}a_k(x)\\[2ex]&=\max_{1\leq
k\leq n}\(\nab V_i(x)\stackrel{\centerdot}{} v_{x_k}\)<0.\ea
$$
This finishes the proof of the lemma. $\Box$ \Vs

\bp\label{p:6.2} Let $W,W'$  be two open attractor neighborhoods of
$\a$. Then there exists a  compact attractor neighborhood $\O$ of
$\a$ such that $\O$ is a strong deformation retract of both $W$ and
$W'$.

Consequently, all  open attractor neighborhoods of $\a$  have the
same homotopy type.

\ep \noindent{\bf Proof.} Let $V,V'$ be smooth Lyapunov functions of
$\a$ on $W$ and $W'$ given by Theorem \ref{t:l2}, respectively. Take
two positive numbers $0<\de<\ve$ sufficiently small such that
$$V'_\de\subset V_{\ve/2}\subset V_\ve\subset W'.$$ We first show that $V_\ve$ and $V'_\de$ are strong deformation retracts
of $W$ and $W'$, respectively.

Indeed, let $S(t)$ be the semiflow on the phase space  $X:=W$
generated by the system: \be\label{e:5.1} x'(t)=-\nab V(x(t)),\Hs
x(t)\in W. \ee Namely, $S(t)x$ is the unique solution of the the
system for each $x\in W$. It is clear that $S(t)$ is well defined on
$W$. Moreover, since $\nab V(x)\ne 0$ outside $\a$, one easily
deduces that $S(t)$ has a global attractor $\cA\subset \a$.

Define a function $t(x)$ on $W$ as: $$
t(x)=\left\{\ba{ll}\sup\{t\geq0|\,\,S([0,t))x\subset W\setminus
V_\ve\},\hs\hs &x\in W\setminus V_\ve;\\[1ex]
0,& x\in V_\ve.\ea\right.$$ As $V_\ve$ is a neighborhood of $\cA$
and $\cA$ attracts $x$, we see that $t(x)$ is finite for each $x\in
W$. Because $\nab V$ and $\pa V_\ve$ is transversal at any point
$x\in \pa V_\ve$, by the basic knowledge on geometric theory of ODEs
we know that $t(x)$ is continuous in $x$.

 Define
$$
H(\sig,x)= S(\sig\, t(x))x,\Hs x\in W.$$ Then $H: [0,1]\X W\ra W$ is
continuous and satisfies: $$H(0,\.)=\mb{id}_{W},\Hs H(1,W)\subset
V_\ve,$$
$$
H(\sig, x)=x,\Hs \A \,(\sig,x)\in [0,1]\X V_\ve.
$$
That is,  $V_\ve$ is a strong deformation retract of $W$.

The same argument applies to show that $V'_\de$ is a strong
deformation retract $W'$.

We claim that $V'_\de$ is a strong deformation retract of $V_\ve$.
For this purpose let $K=V_\ve\-\mb{int}V'_\de$. Then $K\subset
(W\cap W')\-\a$ and is  compact. It follows by Lemma \ref{l:6.1}
that there exists a smooth vector field $\Psi$ defined on $X$ such
that \be\label{e:6.1} \nab V(x)\stackrel{\centerdot}{} \Psi(x)<0,\hs
\nab V'(x)\stackrel{\centerdot}{} \Psi(x)<0 \ee for all  $x\in K$.
 Consider the semiflow $T(t)$
generated by
$$
x'(t)=\Psi(x(t)).
$$
 (\ref{e:6.1}) implies that both $V_\ve$ and $V'_\de$ are
positively invariant with respect to $T(t)$; moreover, the  vector
field $\Psi$ is transversal to both $\pa V_\ve$ and $\pa V'_\de$\,.
Making use of $T(t)$ one can easily construct
 a continuous function $G:
[0,1]\X V_\ve\ra V_\ve$, which is a strong deformation from $V_\ve$
to $V'_\de$ and thus proves our claim. Since the argument is quite
similar as above, we omit the details.

Let $\O=V'_\de$. Recall that $\O$ is a strong deformation retract of
$W'$. In what follows we show that $\O$ is also a strong deformation
retract of $W$, thus finishes the proof of the proposition. For this
purpose we define $\Theta: [0,1]\X W\ra W$ as follows:
$$ \Theta(\sig,x)=\left\{\ba{ll}H(2\sig,x),\hs\hs&0\leq \sig\leq
1/2,\,\,\,x\in W;\\[2ex]
G\(2\sig-1,H(1,x)\),\hs\hs&1/2\leq\sig\leq 1,\,\,\,x\in
W.\ea\right.$$ Clearly $\Theta$ is continuous, and $
\Theta(0,\.)=\mb{id}_{W}$. We observe that  $$
\Theta(1,W)=G(1,H(1,W))\subset G(1,V_\ve)\subset\O.$$ Let
$x\in\O\,\,(\subset V_\ve)$. Then
$$
\Theta(\sig,x)=H(2\sig,x)=x,\Hs \mb{if }\sig\leq 1/2,
$$
and
$$
\Theta(\sig,x)=G(2\sig-1,\,H(1,x))=G(2\sig-1,x)=x,\Hs \mb{if }\sig>
1/2,
$$
Therefore $\Theta$ is a strong deformation from $W$ to $\O$.
 $\Box$

\subsection{Critical groups of Morse sets}

We denote by $H_*$ the usual singular homology theory with
coefficients in a given Abelian group $\mathscr{G}$. Let $\a$ be an
attractor of (\ref{e:1.1})
 with the attraction basin $\W=\W(\a)$, and let  $\cM=\{M_1,\cdots,
 M_l\}$ be a Morse decomposition of $\a$ with Morse filtration
 $$\emp=A_0\subset
A_1\subset\cdots\subset A_l=\a.$$ We first prove the following basic
fact.

\bp\label{t:6.3} Let $W$ and $W'$ be two  open attractor
neighborhoods of $A_k$, and let $U$ and $U'$ be two open attractor
neighborhoods of $A_{k-1}$. Then
$$
H_*(W,U)\cong H_*(W',U').
$$\ep
\noindent{\bf Proof.} Proposition \ref{p:6.2} allows us to pick
compact attractor  neighborhoods $\O$ of $A_k$ and $K$ of $A_{k-1}$
with $K\subset \O$, such that $\O$ (resp. $K$) is a strong
deformation retract of both $W$ and $W'$ (resp. $U$ and $U'$).
Consider the commutative diagram:
$$
\ba{ccccccccc} H_q(K) & \stackrel{i_*}{\longrightarrow}& H_q(\O)&
\stackrel{j_*}{\longrightarrow}&
 H_q\(\O,K\)& \stackrel{\pa}{\longrightarrow}& H_{q-1}(K)& \stackrel{i_*}{\longrightarrow}&
 H_{q-1}(\O)\\[1ex]
\,\,\downarrow{{\tiny i}_* }& &\,\,\downarrow{i_* }&
&\,\,\,\,\downarrow{i_* }& &\,\,\downarrow{i_* }&
&\,\,\downarrow{i_*
}\\[1ex]
H_q(U) & \stackrel{i_*}{\longrightarrow}& H_q(W)
 & \stackrel{j_*}{\longrightarrow}& H_q\(W,U\) &
 \stackrel{\pa}{\longrightarrow}&
H_{q-1}(U)  & \stackrel{i_*}{\longrightarrow}& H_{q-1}(W)
 \ea
$$
The upper and lower rows present the exact homology sequences for
the pairs  $(\O,K)$ and $(W,U)$, respectively. The homomorphisms
 $i_*$'s in the vertical arrows are induced by inclusions.
Since the vertical arrows number 1,2,4 and 5 are isomorphisms, we
immediately
 conclude by the well known ``Five-lemma'' (see \cite{Spa}, Lemma IV.5.11) that $$H_q(W,U)\cong
 H_q(\O,K).$$

 Similarly we also have $$H_q(W',U')\cong
 H_q(\O,K).$$ This completes the  proof of the proposition. $\Box$

\Vs Now let us introduce the concept of {\em critical group } for
Morse sets.

\begin{definition}\label{d:5.1}
The {critical group} $C_*(M_k)$ of Morse set $M_k$ is defined to be
the homology theory  given by
$$
C_q(M_k)=H_q\(W,\,U\),\Hs q=0,1,\cdots, $$ where $W$ and $U$ are
 open attractor neighborhoods of $A_k$ and $A_{k-1}$, respectively.
\end{definition}

\br If $M_k=\emp$, then we necessarily have $A_k=A_{k-1}$. In such a
case both $W$ and $U$ are open attractor neighborhoods of $A_{k-1}$.
Take a $K\subset W\cap U$ so that $K$ is a strong deformation
retract of both $W$ and $U$. Then one finds that
$$C_*(M_k)=H_*(W,U)\cong H_*(K,K)=0.$$

\er

\begin{proposition}\label{p:6.4} Let $V$ be a $C^1$ strict M-L function of
$\cM$ satisfying $(\ref{e:v8})$ and $(\ref{e:3.7})$. Take two real
numbers $a<b$  such that

\begin{enumerate}\item[$(1)$] $M_k$ is the unique Morse set contained in
$V^{-1}([a,b])$; \item[$(2)$]  if $M_k\ne\emp$, then $a<c_k<b$,
where $c_k=V(M_k).$
\end{enumerate}
Then
$$
C_*(M_k)=H_*(V_b,V_a).
$$
\end{proposition}

\noindent{\bf Proof.} Choose a number $\ve>0$ small enough so that
$a+\ve<b-\ve$.  If $M_k\ne\emp$, we also require $a+\ve<c_k<b-\ve$.
Using the semiflow $S(t)$ generated by the gradient system:
$$
x'(t)=-\nab V(x(t)),
$$
it can be shown that $V_{b-\ve}$ is a strong deformation retract of
both $\mb{int}V_b$ and $V_b$,\, and $V_a$ is a strong deformation
retract of
 $\mb{int}V_{a+\ve}$.

Note that $\mb{int}V_b$ and $\mb{int}V_{a+\ve}$ are open attractor
neighborhoods of $A_k$ and $A_{k-1}$, respectively. Hence  by the
definition of critical groups we have
 $$
 C_*(M_k)=H_*\(\mb{int}V_b,\,\mb{int}V_{a+\ve}\).
 $$
On the other hand, using the same argument as in the proof of
Proposition \ref{t:6.3} one easily verifies
$$
H_*\(\mb{int}V_b,\,\mb{int}V_{a+\ve}\)\cong
H_*\(V_{b-\ve},\,V_{a}\)\cong H_*\(V_{b},\,V_{a}\),
$$
and the conclusion of the proposition follows. $\Box$

\br\label{r:6.7} Let $V$, $a, b$ and $M_k$ be the same as in
Proposition \ref{p:6.4}. If \,$V$ has no critical points in $M_k$
$($and hence in $V^{-1}([a,b])$\,$)$, then it can be easily shown
that $V_a$ is a strong deformation retract of $V_b$. It follows that
$ C_*(M_k)=0. $ As a matter of fact, we deduce  in general that
either $ C_*(M_k)=0, $ or $M_k$ contains at least a critical point
of $V$.

A particular but important case is that $M_k$ consists of exactly
one equilibrium $z$ of the system $(\ref{e:1.1})$. In this case if
$C_*(z)\ne0,$, then $z$ is necessarily a critical point of $V$, i.e,
$\nab V(z)=0$. By the basic knowledge in the theory of variational
methods we conclude that \be\label{e:6.5} C_*(z)=H_*(V_b,V_a)\cong
H_*(V_c,V_c\-\{z\}), \ee where $c=c_k=V(z)$; see, for instance,
Chang $\cite{Chang}$. Further by excision of homologies one deduces
that for any neighborhood $U$ of $z$, it holds that \be\label{e:6.6}
C_*(z)\cong H_*(V_c,V_c\-\{z\})\cong H_*\(V_c\cap U,\,V_c\-\{z\}\cap
U\). \ee

Note that if $z$ is asymptotically stable, then $(\ref{e:6.5})$ and
$(\ref{e:6.6})$
imply \be\label{e:6.7}C_*(z)\cong H_*(\{z\})=\left\{\ba{ll}\mathscr{G},\hs\hs &q=0;\\[1ex]
0,&q\geq 1,\ea\right.\ee\er

\vs \noindent{\em Example} 6.1.  Consider the following
 differential inclusion which relates to the generalized
equations governing Chua's circuit \cite{Brown}:
\begin{eqnarray}\label{e:6.3}
\(\ba{ll}\dot{x}_1\\
\dot{x}_2\\
\dot{x}_3\ea\)\in A\(\ba{c}x_1-k\,\mb{Sgn} (x_1)\\
x_2\\
x_3+k\,\mb{Sgn }(x_1)\ea\),\hs A=\(\ba{c}-\alpha(b+1)\\
1\\
0\ea \ba{c}\alpha\\-1\\-\beta\ea \ba{c}0\\1\\0\ea\),
\end{eqnarray}
where  $\mb{Sgn}(x)$ corresponds to  the signal function,
$$
\mb{Sgn}(x)=1\,\,(x>0), \hs
\mb{Sgn}(x)=-1\,\,(x<0),\hs\,\mb{Sgn}(0)=[-1,1].
$$
Taking $ \alpha=-1,$ $ \beta=288,$  $b=-36,$ and $ k=1, $ the system
reads:
\begin{eqnarray}\label{e:6.4}
\(\ba{ll}\dot{x}_1\\
\dot{x}_2\\
\dot{x}_3\ea\)\in A_0\(\ba{c}x_1\\
x_2\\
x_3\ea\)+\(\ba{c}35\mb{Sgn}(x_1)\\
0\\
0\ea\),\hs A_0=\(\ba{c}-35\\
1\\
0\ea \ba{c}-1\\-1\\-288\ea \ba{c}0\\1\\0\ea\).
\end{eqnarray}
Simple computations show that all the eigenvalues of $A_0$ are
negative, so  the system (\ref{e:6.3}) is  dissipative and has a
global attractor $\a$. (\ref{e:6.3}) has three equilibria:
$$
E_1=(-1,0,1),\hs E_2=(1,0,-1),\hs E_3=(0,0,0),
$$
where $E_1$ and $E_2$ are asymptotically stable (therefore each one
is an attractor). Let
$$ A_0=\emptyset,\hs A_1=\{E_1\},\hs A_2=\{E_1,E_2\},\hs A_3=\a.
$$ Then $\{A_k\}$ is an increasing attractor sequence which yields a
Morse decomposition $\cM=\{M_1,M_2,M_3\}$ of $\a$ with
$$
M_1=\{E_1\},\hs M_2=\{E_2\},\hs E_3\in M_3.
$$

For simplicity we take the coefficients group
$\mathcal{G}=\mathbb{Z}$. By (\ref{e:6.7}) one finds that $$
C_q(M_i)=\left\{\ba{ll}\mathbb{Z},\hs\hs &q=0;\\[1ex]
0,&q\geq 1,\ea\right.\Hs i=1,2.
$$
Now let us compute $C_*(M_3)$. Choose open attractor neighborhoods
$U_1$ of $E_1$ and $U_2$ of $E_2$ with $U_1\cap U_2=\emp$. Then
$U=U_1\cup U_2$ is an open attractor neighborhood of $A_2$. We
observe that
$$
H_*(U_i)=C_*(U_i)=H_*(\{E_i\}).
$$
Hence
$$
H_q(U)\cong H_q(U_1)\oplus H_q(U_2)=\left\{\ba{ll}\mathbb{Z}\oplus\mathbb{Z},\hs\hs &q=0;\\[1ex]
0,&q\geq 1.\ea\right.
$$
Let $W=\R^3$. By definition of critical group we have
$$
C_*(M_3)=H_*(W,U).
$$
Using the  exact sequence
$$
\cdots\longrightarrow H_1(W,U) \stackrel{\pa}{\longrightarrow}
H_{0}(U)\stackrel{i_*}{\longrightarrow}H_0(W)\stackrel{j_*}{\longrightarrow}H_0(W,U)\longrightarrow0,
$$
one finds that $$ H_0(W,U)\cong
H_0(W)/\mb{Ker}\,(j_*)=H_0(W)/\mb{Im}\,(i_*). $$ It is easy to see
that $\mb{Im}\,(i_*)=H_0(W)$. Thus we obtain $H_0(W,U)=0$.

To compute $H_1(W,U)$, we consider the long exact sequence of
reduced homologies
$$
\cdots\longrightarrow\~H_q(W)\stackrel{j_*}{\longrightarrow}\~H_q(W,U)
\stackrel{\pa}{\longrightarrow}
\~H_{q-1}(U)\stackrel{i_*}{\longrightarrow}\~H_{q-1}(W)\longrightarrow\cdots.
$$
Noticing that $ \~H_q(W)=\~H_{q-1}(W)=0$ for all $q\geq 0$, one
 deduces that
$$
\~H_q(W,U)\cong \~H_{q-1}(U)=\left\{\ba{ll}\mathbb{Z},\hs\hs &q=1;\\[1ex]
0,&q> 1.\ea\right.
$$
Hence by definition of reduced homologies we conclude that
$$
H_q(W,U)=\~H_q(W,U)=\left\{\ba{ll}\mathbb{Z},\hs\hs &q=1;\\[1ex]
0,&q> 1.\ea\right.
$$
Therefore
$$
C_q(M_3)=\left\{\ba{ll}\mathbb{Z},\hs\hs &q=1;\\[1ex]
0,&q\ne 1.\ea\right.
$$

\subsection{Morse inequalities and Morse equation}
Now let us establish Morse inequalities and Morse equations for
attractors. Let \be\label{e:5.3}
\mathfrak{m}_q=\sum_{k=1}^l\mb{rank}\,C_q\(M_k\),\Hs q=0,1,\cdots.
\ee $\mathfrak{m}_q$ is called the $q-$th {\bf Morse type number }
of $\cM$.

 \bt\label{t:5.2} $(${\bf
Morse inequality and equation}$)$ Let  $\b_q=\b_q(\W):=\mb{\em
rank}\,H_q(\W)$ be  the $q$-th Betti number of the attraction basin
$\W.$ Then the following inequalities and equation hold:
$$
\mathfrak{m}_0\geq \beta_0,$$
$$
\mathfrak{m}_1-\mathfrak{m}_0\geq \beta_1-\beta_0,
$$
$$
\cdots\cdots
$$
$$
\mathfrak{m}_m-\mathfrak{m}_{m-1}+\cdots+(-1)^m\mathfrak{m}_0=\beta_m-\beta_{m-1}+\cdots+(-1)^m\beta_0.
$$

\et \br\label{r:5.1} If we define formal
Poincar$\acute{e}$-polynomials
$$
P_\a(t)=\sum_{q=0}^m \beta_qt^q,\Hs M_\a(t)=\sum_{q=0}^m
\mathfrak{m}_qt^q,
$$
then the Morse inequalities and Morse equation in above can be
reformulated in a very simplified manner:
\be\label{e:5.7}M_\a(t)-P_\a(t)=(1+t)Q_\a(t),\ee where
$Q_\a(t)=\sum_{q=0}^m \gamma_q\,t^q$ is a a formal polynomial with
$\gamma_q$ being nonnegative integers.

\er

To prove Theorem \ref{t:5.2}, we first need to recall  some basic
facts.

A real function $\Phi$ defined on a suitable family $D(\Phi)$ of
pairs of spaces is said to be {\em subadditive}, if $W\subset
Z\subset Y$ implies
$$
\Phi(Y,W)\leq \Phi(Y,Z)+\Phi(Z,W).
$$
 If $\Phi$ is subadditive, then for any  $Y_0\subset Y_1\subset\cdots\subset Y_n$ with  $(Y_k,Y_{k-1})\in D(\Phi)$,
 $$
 \Phi(Y_n,Y_0)\leq \sum_{k=1}^n\Phi(Y_k,Y_{k-1}).$$

For any  pair $(Y,Z)$ of spaces, set
$$
R_q(Y,Z)=\mb{rank}\, H_q(Y,Z)\hs(q\mb{-th {\em Betti number}}).
$$
Define
$$
\Phi_q(Y,Z)=\sum_{j=0}^q(-1)^{q-j}R_j(Y,Z), \Hs
\chi(Y,Z)=\sum_{q=0}^\8(-1)^qR_q(Y,Z).
$$
$\chi(Y,Z)$ is usually called the {\em Euler number} of $(Y,Z)$.

\bl\label{l:5.1} $\cite{Chang, Mil}$  The functions $R_q,\,\Phi_q$
are subadditive, and $\chi$ are additive.\el

\noindent{\bf Proof of Theorem \ref{t:5.2}.} \,We may assume that
all the Morse sets are nonvoid, as the critical group of such a
Morse set is trivial. The following argument is quite standard as in
the case of the classical Morse theory.

Let $V$ be a $C^1$  strict M-L functionof $\cM$, and let
$$c_k=V(M_k),\Hs 1\leq k\leq l.$$ Take
 $a,b\in\R^1$ be such that $$a<c_1<c_2<\cdots<c_l<b.$$ As
 $c_1$ is the minimum of $V$ on $\W$, we have $$\emp=V_a\subset\a\subset V_b.$$
Taking $a_k\in \R^1$ ($k=0,1,\cdots, l$) be such that
$$
a=a_0<c_1<a_1<c_2<a_2<\cdots<c_l<a_l=b,$$ by Lemma \ref{l:5.1} one
immediately deduces that
$$
\sum_{i=1}^l\sum_{j=0}^q(-1)^{q-j}R_j\(V_{a_i},\,V_{a_{i-1}}\)\geq\sum_{j=0}^q(-1)^{q-j}R_j\(V_{a_l},\,V_{a_0}\),
$$
that is, \be\label{e:6.2} \sum_{j=0}^q(-1)^{q-j}\mathfrak{m}_j\geq
\sum_{j=0}^q(-1)^{q-j}R_j(V_b).\ee Noting that $V$ has no critical
point in $V^{-1}([b,+\8))$, making use of the semiflow of the
gradient system
$$
x'(t)=-\nab V(x(t)), \Hs x(t)\in\W, $$ it can be easily shown that
$V_b$ is a strong deformation retract of $\W$. Therefore
$H_*(V_b)=H_*(\W)$, and hence  $R_j(V_b)=\b_q$. This and
(\ref{e:6.2}) justify the Morse inequalities.

To prove  the Morse equation, we  observe that
$$
\chi(V_b,V_a)=\sum_{q=0}^m(-1)^q
R_q(V_b,V_a)=\sum_{q=0}^m(-1)^q\b_q.
$$
The additivity of $\chi$ also yields that
$$\ba{ll}
\chi(V_b,V_a)&=\sum_{i=1}^l\chi\(V_{a_i},\,V_{a_{i-1}}\)\\[2ex]
&=\sum_{i=1}^l\,\sum_{q=0}^m(-1)^q R_q\(V_{a_i},\,V_{a_{i-1}}\)\\[2ex]
&=\sum_{q=0}^m(-1)^q \,\sum_{i=1}^l
R_q\(V_{a_i},\,V_{a_{i-1}}\)=\sum_{q=0}^m(-1)^q \mathfrak{m}_q. \ea
$$
Therefore
$$
\sum_{q=0}^m(-1)^q\b_q=\sum_{q=0}^m(-1)^q \mathfrak{m}_q.
$$
This is precisely what we desired. $\Box$

\br\label{r:5.3} If  $\a$ is the global
 attractor of the flow, then  $\b_q$ is the $q$-th Betti number
  of the phase space  $X=\R^m$. Since $X$ is contractable, we have
$$
H_q(X)=\left\{\ba{ll}\mathcal{G},\hs\hs&q=0;\\[2ex]
0, & q\ne 0.\ea\right.
$$
Taking $\cG=\mathbb{Z}$, one obtains $$ \b_0=1,\hs
\b_q=0\,\,(q>0).$$ Consequently the Morse inequalities and Morse
equation read:
$$
\mathfrak{m}_0\geq 1,$$
$$
\mathfrak{m}_1-\mathfrak{m}_0\geq -1,
$$
$$
\cdots\cdots
$$
$$
\mathfrak{m}_m-\mathfrak{m}_{m-1}+\cdots+(-1)^m\mathfrak{m}_0=(-1)^m.
$$

\er

\noindent{\em Example }6.2.  Consider the global attractor $\a$ of
the system $(\ref{e:6.3})$ with Morse decomposition
$\cM=\{M_1,M_2,M_3\}$. Taking $\cG=\mathbb{Z}$, we have
$$\mathfrak{m}_0=2,\hs \mathfrak{m}_1=1,\hs
\mathfrak{m}_2=\mathfrak{m}_3=0.$$ Therefore  the formal
Poincar$\acute{e}$-polynomials are as follows:
$$
P_\a(t)=\sum_{q=0}^m \beta_qt^q\equiv 1,\Hs M_\a(t)=\sum_{q=0}^m
\mathfrak{m}_qt^q=2+t,
$$
and the equation (\ref{e:5.7}) reads
$$(2+t)-1=(1+t)Q_\a(t),$$ from which one also finds that the  formal polynomial
$Q_\a(t)\equiv 1$.

\Vs\Vs

\end{document}